\documentclass{IET-Conf-Paper}

\usepackage{tikz}
\usetikzlibrary{patterns}
\usepackage{pgfplots}
\usepackage{pifont}
\usepackage{multirow}
\usepackage{makecell}
\usepackage{array}
\usepackage{mhchem} 
\usepackage{threeparttable}
\newcolumntype{P}[1]{>{\centering\arraybackslash}p{#1}}
\pgfplotsset{
	axis lines=middle,
	xlabel=$x$,
	ylabel=$y$,
	no markers,
	samples=51,
	trig format plots=rad,
	x axis in pi/.style={
		scaled x ticks={real:\PI},
		xtick scale label code/.code={},
		xtick distance=pi/#1,
		xticklabel={%
			\pgfmathparse{round(100*\tick)/100}
			\ifdim \pgfmathresult pt = 1pt
			\strut$\pi$%
			\else\ifdim \pgfmathresult pt = -1pt
			\strut$-\pi$%
			\else
			\pgfmathifisint{\pgfmathresult}{%
				\strut$\pgfmathprintnumber[int detect]{\pgfmathresult}\pi$%
			}{%
			\strut$\pgfmathprintnumber{\pgfmathresult}\pi$%
		}
		\fi\fi
	},
	xticklabel style={
		/pgf/number format/.cd,
		frac,
		frac whole=false,
	},
},
}
\begin{document}

\title{HIGH RESOLUTION GENERATION EXPANSION PLANNING CONSIDERING FLEXIBILITY NEEDS: THE CASE OF SWITZERLAND IN 2030}

\author{Elena Raycheva\ad{1}\ad{2}\corr, Jared Garrison\ad{3}, Christian Schaffner\ad{2}, Gabriela Hug\ad{1}}

\address{\add{1}{EEH - Power Systems Laboratory, \add{2}{ESC - Energy Science Center, \add{3}{FEN - Research Center for Energy Networks, ETH Zurich, Zurich, Switzerland}}}
\email{{elena.raycheva}@esc.ethz.ch}}


\begin{abstract}
This paper presents a static generation expansion planning formulation in which operational and investment decisions for a wide range of technologies are co-optimized from a centralized perspective. The location, type and quantity of new generation and storage capacities are provided such that system demand and flexibility requirements are satisfied. Depending on investments in new intermittent renewables (wind, PV), the flexibility requirements are adapted in order to fully capture RES integration costs and ensure normal system operating conditions. To position candidate units, we incorporate DC constraints, nodal demand and production of existing generators as well as imports and exports from other interconnected zones. To show the capabilities of the formulation, high-temporal resolution simulations are conducted on a 162-bus system consisting of the full Swiss transmission grid and aggregated neighboring countries.
\end{abstract}

\maketitle

\section{Introduction}
The integration of large RES capacities has a significant impact on power systems planning as it increases the need for operational flexibility from existing and future units. Consequently, the change in system reserve requirements in response to the expected growth of RES has to be accounted for. The main objective of this work is to present a formulation of the Generation Expansion Planning (GEP) problem which accounts for both present and future flexibility needs and demonstrate its functionality on a real-size power system.\\
\indent The goal of GEP is to determine the optimal investments in new generation and storage technologies over a certain planning horizon, in order to meet load growth and replace decommissioned units. 
A detailed review of GEP in the context of increasing integration of RES is presented in \cite{Generation expansion planning review}. Planning generation expansion in power systems with large shares of RES requires modeling the detailed operational constraints of both existing and candidate technologies providing flexibility such as storages, hydro and thermal generators. Furthermore, reserve constraints have to be included to fully capture the costs of integrating intermittent generation and ensure normal system operating conditions. While \cite{IEEE:Evaluating and planning flexibility in sustainable power systems}-\cite{IEEE:Impact of Operating Reserves} recognize this and handle some of the constraints, they focus primarily on thermal units and do not consider hydro or storages as sources of operational flexibility. In contrast, the present work includes different types of flexibility providers.\\ \indent Further flexibility can be provided through imports and exports from other interconnected zones. This is currently considered to be the most convenient and cheapest way to increase flexibility in regions with reliable grid connectivity \cite{Analyzing operational flexibility of electric power system}. Thereby, it is important to model market-based tie line flow constraints as opposed to the full cross-border line limits to better reflect the realistic ability to import/export. Modeling the grid within the considered zone including its connection to neighboring countries, allows to position candidate units at system nodes of interest and determine favorable locations to alleviate, for example, grid congestions. In \cite{IEEE:Strategic Generation Investment Using a Complementarity Approach} and \cite{Do unit commitment constraints affect generation expansion planning?}, DC constraints are included in the investment model, using simplified transmission systems as test cases. In contrast, our work includes a study of the capability of a real-size power system to evolve and cope with the projected increase in intermittent RES capacities.
We address uncertainties related to renewable production by modeling the reserve provision of both existing and candidate units and capture the increased needs and costs in terms of tertiary system reserve requirements in the optimization problem. This is similar to \cite{Integrating economic and engineering models for future electricity market evaluation: A Swiss case study}, however, in the present work we use nodal dispatch, market-based limits on the cross-border tie lines as well as a unit commitment formulation for the operation of conventional thermal generators, which has been shown to have a significant impact on investment decisions \cite{Do unit commitment constraints affect generation expansion planning?}.\\
\indent An additional novelty of the work lies in the high temporal and spacial resolution of the conducted simulations coupled with detailed modeling of flexibility provision, spanning 1) imports/exports from other zones, 2) operation and 3) reserves. The case study in this paper presents two possible scenarios, namely with and without fixed RES target, for the development of the Swiss generation portfolio in the context of nuclear phase-out and desired increase in RES production. As generation in Switzerland, but also in many other regions, is dominated by hydro capacities, a high temporal resolution for both the validation and generation expansion is needed. While in \cite{IEEE:Evaluating and planning flexibility in sustainable power systems}-\cite{IEEE:Impact of Operating Reserves}, \cite{IEEE:Strategic Generation Investment Using a Complementarity Approach}-\cite{Do unit commitment constraints affect generation expansion planning?} a few representative days or weeks are used, in the present work every other day of the year is simulated with hourly resolution. We validate the short-term operation formulation to ensure agreement with historical data since any discrepancies could impact future investment decisions. The remainder of the paper is organized as follows: the problem formulation is in Section \ref{PF}. The case study and results are in Section \ref{TS} while conclusions are drawn in Section \ref{C}.

\section{Problem Formulation}\label{PF}
In the following problem formulation lowercase letters are used to denote variables and uppercase letters denote parameters. The objective of the generation expansion planning problem is to minimize the sum of the production and investment costs of all existing and candidate generation and storage technologies over the planning horizon $T$, which in the present work is fixed to a single year:
\vspace{-0.15cm}
\begin{align} \label{eq:obj fun}
min 
\underbrace{\sum_{j \in J} \sum_{t \in T} (C_{j}^{prod}p_{j,t} + C_{j}^{su}v_{j,t})}_\text{i)}
+ \underbrace{\sum_{s \in S}\sum_{t \in T} C_{s}^{prod} p^{dis}_{s,t}}_\text{ii)} + \\
\underbrace{\sum_{r \in R} \sum_{t \in T} C_{r}^{prod} p_{r,t}}_\text{iii)} 
+ \underbrace{\sum_{n \in N} \sum_{t \in T} C^{ls} ls_{n,t}}_\text{iv)}
+ \underbrace{\sum_{c \in C} u_c^{inv} I_c}_\text{v)} \nonumber &&
\end{align}
where i) - iii) are the production costs of the set of thermal generators $J$, energy storage systems $S$ and non-dispatchable renewable generators $R$, iv) refers to the load shedding costs at the $N$ system nodes and v) are the investment costs associated with building new units from the set of candidate units $C$. All production costs are assumed to be linear functions of the power generated by the given thermal unit, storage system or renewable generator and the associated operational cost parameter $C^{prod}$. The start-up costs of all thermal generators are expressed as linear functions of the cost parameter $C_{j}^{su}$ and the startup binary variable $v_{j,t}$ and are independent of the time since last shut-down. 
For energy storages, only the operational cost associated with purchasing electricity during charging (pumping) are included. The load shedding at any system node $n$ is the product of the load shedding cost parameter $C^{ls}$ and the load shedding variable $ls_{n,t}$.
In the investment cost formulation, $u_c^{inv}$ denotes the investment decision for each candidate unit $c$. The investment cost $I_c$ is annualized to account for differences in lifetime.  
Expression \eqref{eq:obj fun} is subject to four sets of constraints related to: 1) short-term operation, 2) investments, 3) reserve provision, and 4) transmission system, all of which are described in the following.
\subsection{Short-term Operation} \label{ST Operation}
Short-term operation is modeled by incorporating both the production as well as the reserve provision capabilities of thermal generation, storage and non-dispatchable RES technologies. Table \ref{table reserves} shows which technology types can contribute towards secondary (SCR) and tertiary (TCR) control reserve. Primary reserve is not explicitly modeled as it constitutes less than 10\% of the total hourly reserve quantity and in many western EU countries it does not have to be procured locally \cite{SwissGrid}.
\vspace{-0.4cm}

\begin{table}[htb!]
	\begin{center}
	\caption{Reserve provision per technology}
	\vspace{0.1cm}
	\renewcommand{\arraystretch}{1.1}
	\begin{tabular}{lcc}
		\hline
		\textbf{Technology} & \textbf{SCR $\uparrow \downarrow$} & \textbf{TCR $\uparrow \downarrow$} 
		\\
		\hline
		Thermal: Nuclear/Gas/Coal, etc. & \ding{51} & \ding{51} 
		\\
		Storage: Pumped Hydro/Dam & \ding{51} & \ding{51} 
		\\
		Storage: Battery& \ding{51} & \ding{55} 
		\\
		RES: PV/Wind/Run-of-River & \ding{55} & \ding{55} 
		\\
		\hline
	\end{tabular}
	\label{table reserves}
    \end{center}\vspace*{-18pt}
	\vspace{-0.15cm}
\end{table}
 
\vspace{-0.4cm}
\subsubsection{Conventional Thermal Generators}
The Unit Commitment (UC) constraints of thermal generators are based on the tight and compact formulation in \cite{LUSYM: a Unit Commitment Model formulated as a Mixed-Integer Linear Program} and use three binary variables $u_{j,t}, v_{j,t}, w_{j,t}$, respectively for the on/off status, start up and shut down and one continuous variable $p^{min}_{j,t}$ for the power output above minimum by each unit $j$ in each time period $t$ \cite{LUSYM: a Unit Commitment Model
formulated as a Mixed-Integer Linear Program}. The downward generation constraint is: 
\vspace{-0.2cm}
\begin{equation} \label{eq:pmin}
0 \le p^{min}_{j,t} - (r_{j,t}^{SCR\downarrow} + r_{j,t}^{TCR\downarrow}),\forall j, t
\vspace{-0.2cm}
\end{equation}
where $r_{j,t}^{SCR\downarrow}$ and $r_{j,t}^{TCR\downarrow}$ denote the variables for contribution of each generator towards downward SCR and TCR. The upward generation constraints are given by: {\setlength{\belowdisplayskip}{7pt}
\vspace{-0.2cm}
\begin{flalign} \label{eq:pmax2}
p^{min}_{j,t}  + (r_{j,t}^{SCR\uparrow} + & r_{j,t}^{TCR\uparrow}) \le (P_j^{max} - P_j^{min}) u_{j,t} - \\ & \hspace{-0.7cm}(P_j^{max} - SU_j) v_{j,t} ,\forall t, \forall j \in M^{ut}_j = 1\nonumber &&
\end{flalign}
\vspace{-0.9cm}
\begin{flalign} \label{eq:pmax3}
p^{min}_{j,t} + (r_{j,t}^{SCR\uparrow} + & r_{j,t}^{TCR\uparrow}) \le (P_j^{max} - P_j^{min}) u_{j,t} - \\ & \hspace{-0.7cm}(P_j^{max} - SD_j) w_{j,t+1} ,\forall t, \forall j \in M^{ut}_j = 1\nonumber &&
\end{flalign}
where $P_{j}^{max/min}$ refers to the maximum/minimum power output of the conventional generator $j$ and $SU_j/SD_j$ are its start-up/shut-down capabilities. In case the min. uptime of the generator, $M^{ut}_j$, is two hours or more, a tighter formulation is: 
\vspace{-0.2cm}
\begin{flalign} \label{eq:pmax1}
p^{min}_{j,t} & + (r_{j,t}^{SCR\uparrow} + r_{j,t}^{TCR\uparrow}) \le (P_j^{max} - P_j^{min}) u_{j,t} - \\ & \hspace{-0.7cm}(P_j^{max} - SU_j) v_{j,t} - (P_j^{max} - SD_j) w_{j,t+1} ,\forall t, \forall j \in M^{ut}_j \ge 2\nonumber &&
\end{flalign}
We further use the ramping and min. up/down time constraints from \cite{LUSYM: a Unit Commitment Model formulated as a Mixed-Integer Linear Program}. Due to space constraints we omit their formulation here. The logical relation between the different statuses and the total power generated are given by:
\vspace{-0.3cm}
\begin{equation} \label{eq:logic}
u_{j,t-1} - u_{j,t} + v_{j,t} - w_{j,t} = 0 ,\forall t, \forall j
\vspace{-0.2cm}
\end{equation}
\vspace{-0.4cm}
\begin{equation} \label{eq:pgen}
p_{j,t} = P^{min}_{j} u_{j,t} + p^{min}_{j,t} , \forall t, \forall j
\end{equation}
Planned maintenance is modeled using: 
\vspace{-0.25cm}
\begin{equation} \label{eq:refueling}
u_{j,t} \le S_{j,t}, \forall j \in J^{maint}, \forall t
\vspace{-0.1cm}
\end{equation}
where $S_{j,t}$ is the time series indicating the unit's availability throughout the simulation horizon. The bounds of all previously defined variables are: 
\vspace{-0.1cm}
\begin{equation} \label{eq:within}
p^{min}_{j,t}/r_{j,t}^{SCR/TCR\uparrow\downarrow} \ge 0, \; u_{j,t}/v_{j,t}/w_{j,t} \in [0,1], \; \forall t, \forall j
\vspace{-0.1cm}
\end{equation}
\subsubsection{Energy Storage Units} The operational constraints of each storage unit are modeled with three continuous variables:  $p_{s,t}^{dis}$ and $p_{s,t}^{ch}$ are used for the discharge (turbine) and charge (pump) power and are limited by the maximum discharge/charge power $P_{s,t}^{max,dis/ch}$. The variable energy level $e_{s,t}$ is limited by the energy rating (reservoir energy storage level) $E_s^{max}$ and the final storage level $E_{s,T}$ is set equal to the initial value $E_{s,0}$, i.e. the energy level at the beginning and the end of the year are set equal. It is assumed that each storage system can turn on and produce/consume at maximum discharge/charge power instantaneously:
\vspace{-0.15cm}
\begin{equation} \label{eq:batt1}
0\le p_{s,t}^{dis/ch} \le P_s^{max,dis/ch}, \forall t, \forall s, 
\end{equation}
\vspace{-0.6cm}
\begin{equation} \label{eq:SOC1}
E_s^{min} \le e_{s,t} \le E_s^{max}, \forall t, \forall s  \quad \textrm{and} \quad E_{s,T} = E_{s,T0}, \forall s
\end{equation}
\vspace{-0.55cm}
\begin{equation} \label{eq:SOC2}
e_{s,t} = \underbrace{e_{s,t-1} + \eta_s^{ch} p_{s,t}^{ch} - \frac{p_{s,t}^{dis}}{\eta_s^{dis}}}_\text{$\forall t, \forall s$} \quad +\underbrace{\xi_{s,t}}_\text{$\forall t, \forall s \in S^{hyd}$} 
\end{equation}
\vspace{-0.4cm}
\begin{equation} \label{eq:SOC3}
e_{s,t} \ge 0, \forall t, \forall s
\end{equation}
Eq. \eqref{eq:SOC2} describes the energy content of each storage unit in each hour, taking into account the charging/discharging efficiencies  $\eta_b^{ch}$ and $\eta_s^{dis}$. For hydro power plants, the hourly inflows are included in the term $\xi_{s,t}$. The upward and downward reserve constraints are:
\vspace{-0.15cm} 
\begin{equation} \label{eq:battery reserves1}
r_{s,t}^{SCR\uparrow} + r_{s,t}^{TCR\uparrow}  \le P_{s}^{max,dis}-p_{s,t}^{dis}+p_{s,t}^{ch}, \forall t, \forall s \in S
\vspace{-0.15cm} 
\end{equation}
\begin{equation}  \label{eq:battery reserves2}
r_{s,t}^{SCR\downarrow} + r_{s,t}^{TCR\downarrow} \le P^{max,ch}_{s}-p_{s,t}^{ch}+p_{s,t}^{dis},\forall t,  \forall s \notin S^{dam} \end{equation}
\begin{equation} \label{eq:battery reserves3}
r_{s,t}^{SCR/TCR\uparrow\downarrow} \ge 0, \forall t, \forall s
\end{equation}
where for batteries the variable contribution towards tertiary reserve $r_{s,t}^{TCR\uparrow\downarrow}$ is set to zero. Constraints \eqref{eq:batt1}-\eqref{eq:battery reserves1} are also valid for hydro dams without pumping capabilities with $p_{s,t}^{ch}$ forced to zero. Constraint \eqref{eq:battery reserves1} allows all storage types to provide upward reserve even if they are not producing. Similarly, pumped hydro and batteries can provide downward reserve while staying idle. This assumption is valid as storage units are considered to be infinitely flexible. To ensure that hydro dams do not provide downward reserve when not producing, we add:
\vspace{-0.15cm} 
\begin{equation} \label{eq:dams1}
0 \le p_{s,t}^{dis} - (r_{s,t}^{SCR\downarrow} + r_{s,t}^{TCR\downarrow})  \le P_{s}^{max,dis}, \forall s \in S^{dam}, \forall t 
\vspace{-0.25cm}
\end{equation}
\subsubsection{Non-Dispatchable RES}
Production from solar, wind and run-of-river power plants is modeled via exogenously determined capacity factor profiles, $CF_{r,t}$ multiplied by the unit's maximum installed power $P^{max}_{r}$. We further allow for curtailment of renewable power, i.e.: 
\vspace{-0.2cm}
\begin{equation} \label{eq:RES1}
0\le p_{r,t} \le CF_{r,t} P^{max}_{r},  \forall t, \forall r 
\vspace{-0.2cm}
\end{equation}
\vspace{-0.65cm}
\subsection{Investments} \label{sec:investment constraints}
For thermal generators, the investment decision variable $u_c^{inv}$ from \eqref{eq:obj fun} is binary which corresponds to investments in discrete units. To only dispatch units that have been built, the investment and operational decisions for candidate units are linked:
\vspace{-0.3cm}
\begin{equation} \label{eq:invest1}
u_{c,t} \le u_c^{inv}, \quad u_c^{inv}\in[0,1], \forall t, \forall c \in \mathcal{C}^{thermal}
\vspace{-0.2cm}
\end{equation}
where $u_{c,t}$ is the binary variable for the on/off status of each thermal candidate unit in each time step. Similarly, for storages:
\vspace{-0.2cm}
\begin{equation} \label{eq:invest2}
0 \le p_{c,t}^{dis/ch} \le u_c^{inv}P_{c}^{max,dis/ch},  u_c^{inv}\in[0,1], \forall t, \forall c \in \mathcal{C}^{storage}
\vspace{-0.1cm}
\end{equation}
To only allow reserve provision by units that are built, we use the same constraints as for already existing units, but multiply $P_{s}^{max,dis/ch}$ in \eqref{eq:battery reserves1}-\eqref{eq:battery reserves2} by $u_{c}^{inv}$. For non-dispatchable RES candidate generators, the investment decision variable $u_c^{inv}$ is continuous and corresponds to the built capacity at the candidate location with capacity factor $CF_{c,t}$ and maximum allowable investment capacity $P^{inv,max}$:   
\vspace{-0.15cm}
\begin{equation} \label{eq:invest3}
0 \le p_{c,t} \le u_c^{inv} CF_{c,t}, \quad 0 \le u_c^{inv} \le P_{c}^{inv,max}, \forall t, \forall c \in \mathcal{C}^{RES},
\end{equation} 
\vspace{-0.65cm}
\subsection{Reserve Provision} \label{sec:reserves section}
The formulation of the reserve constraints in Section \ref{ST Operation} allows for non-symmetric reserve provision by each generator/storage unit which is consistent with efforts to reduce market barriers for smaller bidders who might be unable to offer symmetrical power bids \cite{SwissGrid}. The reserves provided by the units have to satisfy the system-wide demand for up/down balancing capacity in each time period:
\vspace{-0.2cm}
\begin{equation} \label{eq:system reserve 0}
\sum_{j \in J} r_{j,t}^{TCR\uparrow} + \sum_{s \in S^{hydro}} r_{s,t}^{TCR\uparrow} \ge TCR_{t}^{\uparrow, sys}  + r^{TCR \uparrow, RES}, \forall t
\end{equation}
\begin{equation} \label{eq:system reserve 1}
\sum_{j \in J} r_{j,t}^{TCR\downarrow} + \sum_{s \in S^{hydro}} r_{s,t}^{TCR\downarrow} \ge TCR_{t}^{\downarrow, sys} + r^{TCR \downarrow, RES},\forall t
\vspace{-0.1cm}
\end{equation}
where $TCR_{t}^{\uparrow, sys}$ is the upward tertiary system reserve quantity required by the Transmission System Operator (TSO) at time step $t$. Depending on the investments in wind and solar PV capacities, an additional tertiary reserve quantity $r^{TCR \uparrow \downarrow, RES}$ is added to ensure that there is enough system flexibility to compensate uncertainties in RES production:
\vspace{-0.2cm}
\begin{equation} \label{eq:system reserve 2}
r^{TCR \uparrow \downarrow, RES} = A_{wind}^{\uparrow \downarrow} \sum_{c \in C_{wind}^{RES}} u_c^{inv} + A_{pv}^{\uparrow \downarrow} \sum_{c \in C_{pv}^{RES}} u_c^{inv}
\vspace{-0.1cm}
\end{equation}
where $A_{wind/pv}$ is an empirically derived coefficient calculated following the methodology in \cite{Integrating economic and engineering models for future electricity market evaluation: A Swiss case study} where short-term wind and PV forecast methods were used to quantify the additional reserves needed. The constraints for provision of secondary reserve are identical to \eqref{eq:system reserve 0}-\eqref{eq:system reserve 1} without the additional terms from \eqref{eq:system reserve 2}. Similar to~\cite{ZimaReserves}, this formulation assumes that the variability in RES generation is accounted for in the tertiary reserve requirement.
\vspace{-0.15cm} 
\subsection{Transmission System} \label{sec:grid constraints}
The following equation models the active power balance at each bus node $n \in N$ where $P_{n,t}^{D}$ is the nodal demand, $ls_{n,t}$ refers to the load shedding variable, and the remaining terms correspond to the power output of each generator or storage:
\vspace{-0.25cm}
\begin{flalign} \label{eq:nodal balance1}
p_{n,t} & = P_{n,t}^{D}-ls_{n,t}+\sum_{s\in S_{n,t}}p_{s,t}^{ch}-\\&\sum_{j\in J_{n,t}} p_{j,t}-\sum_{s\in S_{n,t}}p_{s,t}^{dis}-\sum_{r\in R_{n,t}}p_{r,t}, \quad \forall n,\forall t \nonumber &&
\vspace{-0.4cm}
\end{flalign}
The nodal active power $p_{n,t}$ is the sum of the active power flows of all lines $l \in L$ connected to $n$ as given in:
\vspace{-0.2cm}
\begin{equation} \label{eq:nodal balance3}
p_{n,t} = \sum_{i\in l(n,i)} p_{l(n,i),t}, \forall t, \forall n
\vspace{-0.2cm}
\end{equation}
and the active power flow $p_{l}$ of a single line is:
\vspace{-0.25cm}
\begin{equation} \label{eq:nodal injection1}
p_{l(n,i),t} = B_{l}(\delta_{n,t} - \delta_{i,t}), \forall t
\vspace{-0.1cm}
\end{equation}
\begin{equation} \label{eq:thermal limit}
-P_{l}^{max} \le p_{l(n,i),t} \le P_{l}^{max},\forall t,\forall l(n,i), 
\end{equation}
where $B_{l}$ is the admittance, $\delta_n$, $\delta_i$ are the voltage angles at the start and end nodes and $P_{l}^{max}$ is the thermal limit of the line. 
Load shedding is allowed at each demand bus:
\vspace{-0.2cm}
\begin{equation} \label{eq:nodal balance2}
0 \le ls_{n,t} \le P_{n,t}^{D},\forall t, \forall n 
\vspace{-0.2cm}
\end{equation}
\vspace{-0.6cm}
\section{Case Study} \label{TS}
\subsection{Test System}
In this section, the proposed formulation is applied to an interconnected system consisting of the detailed Swiss transmission network with aggregated neighboring countries. This model in total comprises 263 transmission lines, 162 substations, 21 transformers and 355 existing generators with total installed capacity of 460 GW in 2015. As Switzerland's generation portfolio is heavily dominated by hydro capacities, capturing their operational behavior is salient to any model attempting to replicate historical or predict future production. To speed up computations, while maintaining very high temporal resolution (necessary due to short-term fluctuations in river flows, wind and solar generation) and chronological accuracy (necessary due to the presence of seasonal storages), every other day of the year is simulated with hourly resolution. Thus, the change in demand behavior between weekdays and weekends during each week is always captured.  \\
\indent Fig.\ref{DaysCompression} shows how hydro storage levels are approximated for the days which are not simulated. This form of compression is only used for pumped and dam hydro power plants and not for battery storages as it is assumed the latter operate on shorter cycles (less than a day). Our approach relies on the assumption of day-to-day similarity in operation of both pump and dam power plants. This is valid for dams as they operate on a seasonal cycle as well as for pumped power plants which, depending on their reservoir capacity, operate on a daily to weekly cycle. By adapting \eqref{eq:SOC1}-\eqref{eq:SOC2}, the pumping/turbining across two days is aggregated into the time during which the storage charges/discharges in a single day, which means that the modeled fluctuations in storage level would have double the amplitude. This doubling is not relevant for seasonal storages, but is relevant for those that operate on a daily cycle. Therefore, the initial/minimum/maximum reservoir levels of daily pumped storages are doubled. \\
\indent To establish the investments and UC schedule, the MILP formulation is implemented in Pyomo \cite{Pyomo} and solved with Gurobi \cite{Gurobi}. The considered investment costs span the entire year instead of only every second day, therefore we double the operating costs in \eqref{eq:obj fun}. As it is possible to get non-unique solutions for the hourly operation of hydro storages, stemming from the aggregated modeling of the surrounding countries and simplified production costs, we fix the investment and binary UC decisions and re-solve the linear dispatch problem while also including a negligible price incentive for keeping more water in the storages as a security measure. In this way, we are able to choose a specific storage curve out of the ones that all lead to the same objective function value.
\vspace{-0.4cm}
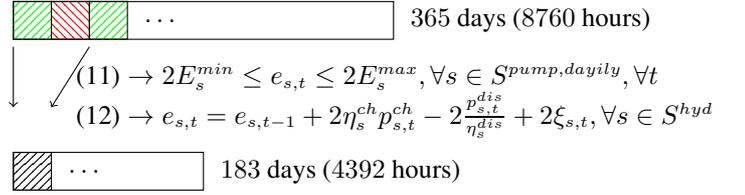
\begin{figure}[!htb]
	\centering
	\begin{tikzpicture}
	\draw (0,4) rectangle (5,3.5);
	\node [right] at (5.1,3.75) {$365$ days ($8760$ hours)};
	
	\draw [pattern color=green,pattern=north east lines] (0,4) rectangle (0.5,3.5);
	\draw [pattern color=red,pattern=north west lines] (0.5,4) rectangle (1.0,3.5);
	\draw [pattern color=green,pattern=north east lines] (1.0,4) rectangle (1.5,3.5);
	\node [right] at (1.6,3.75) {$\ldots$};
	
	\draw [->] (0,3.4) -- (0,2.6);
	\draw [->] (1.0,3.4) -- (0.5,2.6);
	
	\draw (0.0,1.5) rectangle (2.5,2);
	\draw [pattern color=black,pattern=north east lines] (0,2) rectangle (0.5,1.5);
	\node [right] at (0.6,1.75) {$\ldots$};
	\node [right] at (2.6,1.75) {$183$ days ($4392$ hours)};
	
	\node [right] at (0.7,3.0) 
	{$\eqref{eq:SOC1} \rightarrow 2E_s^{min} \le e_{s,t} \le 2E_s^{max}, \forall s \in S^{pump,dayily}, \forall t$};
	\node [right] at (0.7,2.5) 
	{$\eqref{eq:SOC2} \rightarrow e_{s,t} = e_{s,t-1} + 2 \eta_s^{ch}  p_{s,t}^{ch} - 2 \frac{p_{s,t}^{dis}}{\eta_s^{dis}} + 2\xi_{s,t}, \forall s \in S^{hyd}$};
	\end{tikzpicture}
	\vspace{-0.6cm}
	\caption{Days compression for simulation speed-up}\label{DaysCompression}
	\vspace{-0.2cm}
\end{figure}
\vspace{-0.4cm}
\subsection{Validation of Operational Model}
A validation for the year 2015 is performed. The generators in the surrounding countries are aggregated to one unit per technology type. The total capacity per country in 2015 and the operational constraints used are indicated in Table \ref{2015 gen capacities}. The hourly load data and reserve requirements for CH are taken from \cite{swissgrid2}. Load data for the neighbors are from \cite{entsoe}. The cross-border flows between the surrounding and all other countries which are not currently modeled (e.g. DE-DK, etc.) are fixed to the 2015 values from \cite{entsoe}. Solar irradiation and wind time series are from \cite{Integrating economic and engineering models for future electricity market evaluation: A Swiss case study}. The initial/final storage levels for 2015 and hydro inflows are from \cite{BFE}.
\vspace{-0.7cm}\\
\begin{table}[!htbp]
	\centering
	\caption{Generation capacities per country (2015)}
	\vspace{0.1cm}
	\label{2015 gen capacities}
	\renewcommand{\arraystretch}{1.1}
	\begin{threeparttable}
	\begin{tabular}{P{23mm}|ccc|c}
		\hline
		\textbf{Country} & \textbf{Detail} & \textbf{Gens.} & \textbf{Cap. [GW]} & \textbf{Constr.}%
		\\
		\hline
		Austria (AT) & \ding{55} & 9 & 22 & 
		\multirow{4}{*}{\shortstack{Conv.\tnote{*}\\\eqref{eq:batt1}-\eqref{eq:SOC3}\\\eqref{eq:RES1}}}
		\\
		\cline{2-4}
		Germany (DE) & \ding{55} & 12 & 187 &
		\\
		\cline{2-4}
		France (FR) & \ding{55} & 11 & 110 &
		\\
		\cline{2-4}
		Italy (IT) & \ding{55} & 10 & 122 &
		\\
		\cline{1-5}
		Switz. (CH) & \ding{51} & 313 &19 & \eqref{eq:pmin}-\eqref{eq:system reserve 2}
		\\
		\hline
	\end{tabular}
	\begin{tablenotes}
    \item[*] For these conventional units we use linear min/max constraints for $p_{j,t}$, i.e. $ P_{j}^{min} \le p_{j,t} \le P_{j}^{max}$, without modeling ramping and minimum up/down times.
	\end{tablenotes}
\end{threeparttable}
\label{table2}
\end{table}

\vspace{-0.8cm}
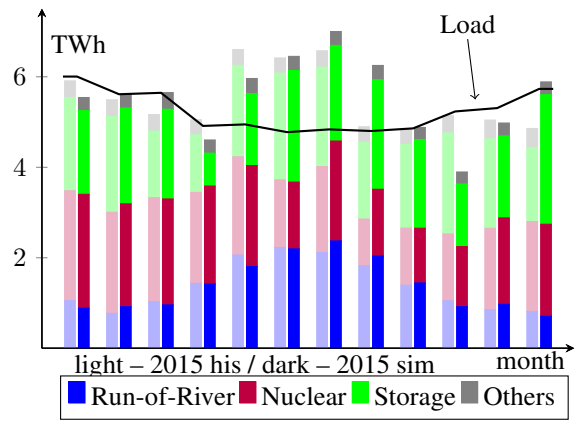
\begin{figure}[!htbp]
	\centering
\begin{tikzpicture}
\node [below] at (2.8,0.04) {light -- 2015 his / dark -- 2015 sim};
\begin{axis}[
ybar stacked,
bar width=4pt,
enlargelimits=0,
ylabel={},
xlabel={},
ytick=\empty,
xtick=\empty,
ymin=0,
xmin=0,
xmax=38,
ymax=7.5,
  scale only axis=true,
  width=7cm,
  height=4.5cm
]
\addplot +[blue!30!white] coordinates {(2,1.049) (5,0.768) (8,1.026) (11,1.433) (14,2.058) (17,2.225) (20,2.115) (23,1.822) (26,1.392) (29,1.046) (32,0.85) (35,0.811) };
\addplot +[purple!30!white] coordinates {(2,2.43) (5,2.235) (8,2.299) (11,2.01) (14,2.172) (17,1.494) (20,1.899) (23,1.032) (26,1.261) (29,1.478) (32,1.8) (35,1.985) };
\addplot +[green!30!white] coordinates {(2,2.061) (5,2.116) (8,1.474) (11,1.271) (14,2.02) (17,2.36) (20,2.199) (23,1.695) (26,1.849) (29,2.232) (32,1.975) (35,1.639) };
\addplot +[gray!30!white] coordinates {(2,0.367) (5,0.37) (8,0.363) (11,0.334) (14,0.349) (17,0.335) (20,0.359) (23,0.344) (26,0.331) (29,0.387) (32,0.416) (35,0.421) };

\end{axis}
\begin{axis}[
ybar stacked,
bar width=4pt,
enlargelimits=0,
legend style={
    legend columns=4,
    at={(0.99,-0.08)},
},
ylabel={TWh},
x label style={at={((1.0,-0.09))}},
y label style={at={(0.0,0.95)}},
xlabel={month},
xtick=\empty,
ymin=0,
xmin=0,
xmax=38,
ymax=7.5,
  scale only axis=true,
  width=7cm,
  height=4.5cm
]
\addplot +[blue] coordinates {(3,0.886029) (6,0.918214) (9,0.958766) (12,1.421518) (15,1.807452) (18,2.196378) (21,2.372435) (24,2.039136) (27,1.444603) (30,0.919237) (33,0.971389) (36,0.706453) };
\addplot +[purple] coordinates {(3,2.516947) (6,2.273371) (9,2.341455) (12,2.163562) (15,2.230925) (18,1.476297) (21,2.207879) (24,1.474916) (27,1.208587) (30,1.327717) (33,1.910157) (36,2.034856) };
\addplot +[green] coordinates {(3,1.848498) (6,2.119242) (9,1.982372) (12,0.727107) (15,1.589511) (18,2.467296) (21,2.11317) (24,2.422647) (27,1.96094) (30,1.384446) (33,1.812448) (36,2.862188) };
\addplot +[gray] coordinates {(3,0.290443) (6,0.298705) (9,0.368073) (12,0.291338) (15,0.334458) (18,0.311253) (21,0.306813) (24,0.313924) (27,0.265869) (30,0.264554) (33,0.283984) (36,0.283004) };

\legend{Run-of-River,Nuclear,Storage,Others}
\end{axis}
\begin{axis}[
enlargelimits=0,
ylabel={},
xlabel={},
ytick=\empty,
xtick=\empty,
ymin=0,
xmin=0,
xmax=38,
ymax=7.5,
  scale only axis=true,
  width=7cm,
  height=4.5cm
]
\addplot +[thick, black] coordinates {(1.5,6.00248) (2.5,6.00248) (5.5,5.61162) (8.5,5.64309) (11.5,4.911502) (14.5,4.944916) (17.5,4.774052) (20.5,4.834574) (23.5,4.799468) (26.5,4.857219) (29.5,5.231399) (32.5,5.305808) (35.5,5.729983) (36.5,5.729983) };

\end{axis}

\node at (5.6,4.3) {Load};
\draw[->] (5.6, 4.1) -- (5.7, 3.3);
\end{tikzpicture}
	\vspace{-0.2cm}
	\caption{Monthly production per technology type in CH (2015)}\label{Monthly Production}
	\vspace{-0.4cm}
\end{figure}\\
\indent Fig. \ref{Monthly Production} compares the 2015 monthly simulated production in Switzerland to the historical values from \cite{BFE}. The simulation results show quite good agreement with the historical production. The most notable difference is the generation of hydro storages in August, October and December. In the simulation, more water is stored in the period April-July and is then turbined in August. This can also be seen in Fig. \ref{Monthly Storage} which shows the water level over the year. Water levels peak at the end of September. Again more water is stored in October and November and used in December. Historically, the dispatch is more evenly distributed. A potential explanation is that the simulation has perfect foresight about demand, inflows and renewable production, which power plant operators do not, therefore they enter long-term contracts to hedge their production. Storing/producing significantly more in any particular month exposes the hydro generators to risks.\\
\vspace{-0.5cm}
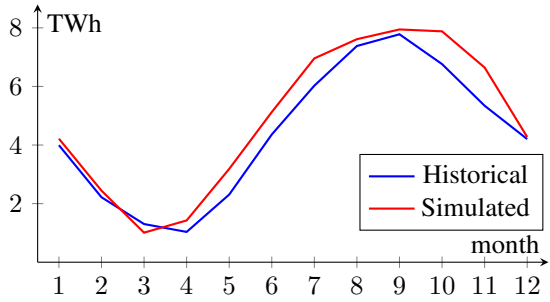
\begin{figure}[!htbp]
	\centering
\begin{tikzpicture}
\begin{axis}[
  enlargelimits=0,
  legend style={
    legend columns=1,
    at={(0.99,0.42)},
  },
  ylabel={TWh},
  xlabel={month},
  xtick={1,2,3,4,5,6,7,8,9,10,11,12},
  ymin=0,
  xmin=0.5,
  xmax=12.5,
  ymax=8.8,
  height=5cm,
  width=8.3cm
]

\addplot +[thick, sharp plot] coordinates {(1,3.995) (2,2.217) (3,1.303) (4,1.034) (5,2.309) (6,4.36) (7,6.027) (8,7.381) (9,7.781) (10,6.763) (11,5.342) (12,4.204) };
\addplot +[thick, sharp plot] coordinates {(1,4.21572715472) (2,2.44301413663) (3,1.00871693838) (4,1.42585878057) (5,3.18537076166) (6,5.1259970097) (7,6.95855561736) (8,7.61366975697) (9,7.94610507168) (10,7.88198473949) (11,6.64620185593) (12,4.27558155812) };

\legend{Historical,Simulated}
\end{axis}
\end{tikzpicture}
    \vspace{-0.2cm}
	\caption{End-of-month storage levels in CH (2015)}\label{Monthly Storage}
	\vspace{-0.4cm}
\end{figure}\\
\indent As future investments in Switzerland could be influenced by the export/import with its neighbours, it is also important to correctly reproduce cross-border flows. Table \ref{exports imports} shows the net Swiss cross-border exchange in 2015. The model is able to capture the most important trend - Switzerland is a transit country for flows from Germany, France and Austria to Italy. \\
\vspace{-0.7cm}
\begin{table}[htb!]
	\centering
	\caption{Net CH cross-border exchange (2015)}
	\vspace{0.1cm}
	\label{exports imports}
	\renewcommand{\arraystretch}{1.1}
	\begin{tabular}{ccc}
		\hline
		\textbf{Net Export (From - To)} & \textbf{Hist. [TWh]} & \textbf{Sim. [TWh]} 
		\\
		\hline
		AT - CH  & 6.7 & 4.9 
		\\
		DE - CH & 13.1 & 9.0
		\\
		FR - CH & 5.3 & 4.3
		\\
		CH - IT & 25.4 & 21.0
		\\
		\hline
	\end{tabular}
\end{table}
\vspace{-0.4cm}
\\Possible reasons for the deviations in absolute magnitude of net cross-border exchange include: the generators and transmission system of the surrounding countries are modeled in aggregate and the assumed variable operation costs are static relative to the actual ranges of operating costs. Regarding the grid topology, all CH cross-border lines go to a single node within the neighboring country which means that in each hour there is a single direction of cross-border exchange, which is not the case in reality. Overall, the presented validation results reflect the historical values well and give confidence in the model. 
\subsection{2030 Generation Expansion Planning}
The expansion planning is established for Switzerland for the target year 2030. Following a 50-year decommissioning plan, only 36\% (1220 MW) of the 2015 installed nuclear capacity in the country will remain operational in 2030. We present two different scenarios: 1) business-as-usual (BaU) and 2) renewable target (RES target). In 2) a production target of 9 TWh from non-hydro renewable generators (including existing biomass, PV and wind) is imposed in Switzerland. In both scenarios, 65 candidate units with varying sizes and cost parameters from \cite{PSI} and own calculations, summarized in Table \ref{2030 candidate units}, are placed at system nodes of interest. While the costs of biomass reflect on-going waste incineration subsidies which are expected to continue in the future, we assume no subsidies for PV and wind. All planned hydro power and transmission system upgrades in the period 2016-2025 are included. Swiss demand and fuel cost projections for 2030 are from \cite{system adequacy} and \cite{primes}. Hydro inflows are set to the 2015 values and the production profiles for PV and wind candidates are from \cite{Integrating economic and engineering models for future electricity market evaluation: A Swiss case study}.\\
\vspace{-0.8cm}
\begin{table}[!htbp]
	\centering
	\caption{Cost parameters of candidate units in CH (2030)}
	\vspace{0.1cm}
	\label{2030 candidate units}
	\renewcommand{\arraystretch}{1.05}
	\begin{tabular}{P{1.1cm}P{2.1cm}P{1.7cm}P{1.1cm}P{0.7cm}}
		\hline
		\textbf{Tech.} & \textbf{Invest. Cost [kEUR/MW/a]}  & \textbf{Var. Cost [EUR/MWh]} & \textbf{Cap. [MW]} & \textbf{Units}
		\\
		\hline
		Gas CC & 84  & 85 & 4200 &28
		\\
		Gas SC & 54  & 131.5 & 600 &14
		\\
		Biomass & 125  & 1 & 240 &12
		\\
		Wind & 206  & 2.5 & 1510 &7
		\\
		PV & 106 & 2.1 & 10000 &4
		\\
		\hline
	\end{tabular}
\end{table}

\vspace{-0.4cm}
\\
\indent The demand, generators and fuel costs in the surrounding countries are adapted to reflect the 2030 projections from \cite{primes}. The 2015 wind and solar production profiles of the neighbors are scaled to match the projected 2030 totals from \cite{primes} and the cross-border flows with all other countries are fixed to the values for 2015 \cite{entsoe}. Table \ref{investments 2030} summarizes the investments made under the two considered scenarios. Even without a RES target, all biomass candidate power plants are built. Given their low costs and the decreased nuclear production, it is more economically viable to have new generators produce locally than to solely import. To satisfy the target, a total of 240 MW biomass and 3254 MW PV is invested in (the remaining 3.36 TWh are produced by existing generators: 2.1 TWh (biomass), 1.1 TWh (PV) and 0.16 TWh (wind)). As a result of the increased intermittent RES generation in the second scenario, the total tertiary system reserve requirements in each hour increase by 26 MW (up) and 28 MW (down). Even with this increase, investments in new dispatchable units were not needed. 
\vspace{-0.4cm}
\begin{table}[!htbp]
	\centering
	\caption{New investments in CH (2030)}
	\vspace{0.1cm}
	\label{investments 2030}
	\renewcommand{\arraystretch}{1.05}
	\begin{tabular}{P{0.5cm}P{1.5cm}P{1.0cm}P{1.0cm}P{1.2cm}P{1.2cm}}
		\hline
		\textbf{Scen.} & \textbf{Techn.} & \textbf{Built [MW]} &  \textbf{Gen. [TWh]} & \textbf{+ TCR$\uparrow$ [MW]} & \textbf{+ TCR$\downarrow$ [MW]}
		\\
		\hline
		BaU  & Biomass & \makecell{240} & \makecell{2.0} & \makecell{\ding{55}} & \makecell{\ding{55}}
		\\
		\hline
		RES & \makecell{Biomass\\PV} & \makecell{240\\3254} & \makecell{2.0\\3.64} & \makecell{\ding{55}\\26} & \makecell{\ding{55}\\28}
		\\
		\hline
	\end{tabular}
\end{table} 
\vspace{-0.4cm}\\
\indent Table \ref{el prices} shows the percentage decrease (-\%) in net generation and increase (+\%) in average annual electricity price in 2030 compared to 2015. The reasons for the price increase in 2030 are twofold: 1) less domestic generation and 2) projected increase in \ce{CO2} and fuel costs. Since Switzerland is a price taker during the majority of the year, the generation costs of the conventional units in the surrounding countries have a profound impact on Swiss electricity prices.\\
\begin{table}[htb!]
	\centering
	\caption{Change in net gen. and av. el. price in CH (2030)}
	\vspace{0.1cm}
	\label{el prices}
	\renewcommand{\arraystretch}{1.1}
	\begin{tabular}{ccc}
		\hline
		\textbf{Scen.} & \textbf{Tot. net gen. [\% 2015]} & \textbf{Av. el. price [\% 2015]}%
		\\
		\hline
		BaU  & -19\% & +51\% 
		\\
		RES & -14\% & +47\%
		\\
		\hline
	\end{tabular}
\vspace{-0.2cm}
\end{table}
\vspace{-0.4cm}
\begin{figure}[!htbp]
	\centering
\begin{tikzpicture}
\node [below] at (3.3,-1.1) {light -- 2015 / normal - 2030 BaU / striped -- 2030 RES};
\begin{axis}[
  ybar stacked,
  bar width=3pt,
  enlargelimits=0,
  ylabel={},
  xlabel={},
  ytick=\empty,
  xtick=\empty,
  ymin=0,
  xmin=0,
  xmax=50,
  ymax=7,
  scale only axis=true,
  width=7cm,
  height=4.5cm
]
\addplot +[blue!30!white] coordinates {(1,0.886029) (5,0.918214) (9,0.958766) (13,1.421518) (17,1.807452) (21,2.196378) (25,2.372435) (29,2.039136) (33,1.444603) (37,0.919237) (41,0.971389) (45,0.706453) };
\addplot +[purple!30!white] coordinates {(1,2.516947) (5,2.273371) (9,2.341455) (13,2.163562) (17,2.230925) (21,1.476297) (25,2.207879) (29,1.474916) (33,1.208587) (37,1.327717) (41,1.910157) (45,2.034856) };
\addplot +[green!30!white] coordinates {(1,1.848498) (5,2.119242) (9,1.982372) (13,0.727107) (17,1.589511) (21,2.467296) (25,2.11317) (29,2.422647) (33,1.96094) (37,1.384446) (41,1.812448) (45,2.862188) };
\addplot +[orange!30!white] coordinates {(1,0.16999) (5,0.153537) (9,0.16999) (13,0.164498) (17,0.16999) (21,0.164498) (25,0.16999) (29,0.16999) (33,0.164498) (37,0.16999) (41,0.164498) (45,0.16999) };
\addplot +[yellow!30!white] coordinates {(1,0.038637) (5,0.065234) (9,0.112132) (13,0.117789) (17,0.154808) (21,0.136414) (25,0.127447) (29,0.136227) (33,0.093281) (37,0.070284) (41,0.045783) (45,0.031855) };
\addplot +[gray!30!white] coordinates {(1,0.009713) (5,0.011453) (9,0.011077) (13,0.009051) (17,0.00966) (21,0.010341) (25,0.009376) (29,0.007707) (33,0.00809) (37,0.00763) (41,0.006809) (45,0.007939) };
\addplot +[teal!30!white] coordinates {(1,0.072103) (5,0.068481) (9,0.074874) (13,0.0) (17,0.0) (21,0.0) (25,0.0) (29,0.0) (33,0.0) (37,0.01665) (41,0.066894) (45,0.07322) };

\end{axis}
\begin{axis}[
  ybar stacked,
  bar width=3pt,
  enlargelimits=0,
  ylabel={},
  xlabel={},
  ytick=\empty,
  xtick=\empty,
  ymin=0,
  xmin=0,
  xmax=50,
  ymax=7,
  scale only axis=true,
  width=7cm,
  height=4.5cm
]
\addplot +[blue,postaction={pattern=north east lines}] coordinates {(3,0.886029) (7,0.918217) (11,0.958766) (15,1.421518) (19,1.806685) (23,2.195713) (27,2.372115) (31,2.039136) (35,1.444603) (39,0.919237) (43,0.97139) (47,0.706454) };
\addplot +[purple,postaction={pattern=north east lines}] coordinates {(3,0.907679) (7,0.819839) (11,0.907679) (15,0.878399) (19,0.907679) (23,0.086863) (27,0.877423) (31,0.907679) (35,0.878399) (39,0.907679) (43,0.878399) (47,0.907679) };
\addplot +[green,postaction={pattern=north east lines}] coordinates {(3,2.096769) (7,2.214424) (11,1.626076) (15,1.469319) (19,1.410887) (23,1.389602) (27,2.065171) (31,2.368246) (35,2.245146) (39,1.718197) (43,1.373692) (47,2.070982) };
\addplot +[orange,postaction={pattern=north east lines}] coordinates {(3,0.348556) (7,0.314823) (11,0.348556) (15,0.33732) (19,0.348556) (23,0.33732) (27,0.348556) (31,0.348556) (35,0.33732) (39,0.348556) (43,0.33732) (47,0.348556) };
\addplot +[yellow,postaction={pattern=north east lines}] coordinates {(3,0.163112) (7,0.275371) (11,0.473306) (15,0.497191) (19,0.653428) (23,0.575799) (27,0.537944) (31,0.575014) (35,0.393751) (39,0.296687) (43,0.193272) (47,0.134496) };
\addplot +[gray,postaction={pattern=north east lines}] coordinates {(3,0.009713) (7,0.011453) (11,0.011077) (15,0.009052) (19,0.00966) (23,0.010341) (27,0.009376) (31,0.007707) (35,0.00809) (39,0.00763) (43,0.006809) (47,0.007939) };
\addplot +[teal,postaction={pattern=north east lines}] coordinates {(3,0.0) (7,0.0) (11,0.0) (15,0.0) (19,0.0) (23,0.0) (27,0.0) (31,0.0) (35,0.0) (39,0.0) (43,0.0) (47,0.0) };

\end{axis}
\begin{axis}[
  ybar stacked,
  bar width=3pt,
  enlargelimits=0,
  legend style={
    legend columns=4,
    at={(0.835,-0.01)},
  },
  ylabel={TWh},
  x label style={at={(1.0,-0.09)}},
  xlabel={month},
  xtick=\empty,
  ymin=0,
  xmin=0,
  xmax=50,
  ymax=7,
  scale only axis=true,
  width=7cm,
  height=4.5cm
]
\addplot +[blue] coordinates {(2,0.886029) (6,0.918214) (10,0.958766) (14,1.421518) (18,1.807452) (22,2.196378) (26,2.372435) (30,2.039136) (34,1.444603) (38,0.919237) (42,0.971389) (46,0.706453) };
\addplot +[purple] coordinates {(2,0.907679) (6,0.819839) (10,0.907679) (14,0.878399) (18,0.907679) (22,0.086863) (26,0.877423) (30,0.907679) (34,0.878399) (38,0.907679) (42,0.878399) (46,0.907679) };
\addplot +[green] coordinates {(2,2.021908) (6,2.080215) (10,1.69883) (14,1.463687) (18,1.634444) (22,1.94911) (26,1.649981) (30,2.057118) (34,2.317023) (38,2.270519) (42,1.170821) (46,1.482354) };
\addplot +[orange] coordinates {(2,0.348538) (6,0.314805) (10,0.348538) (14,0.337286) (18,0.348538) (22,0.337286) (26,0.348538) (30,0.348538) (34,0.337286) (38,0.348538) (42,0.337286) (46,0.348538) };
\addplot +[yellow] coordinates {(2,0.038637) (6,0.065234) (10,0.112132) (14,0.117789) (18,0.154808) (22,0.136414) (26,0.127447) (30,0.136227) (34,0.093281) (38,0.070284) (42,0.045783) (46,0.031855) };
\addplot +[gray] coordinates {(2,0.009713) (6,0.011453) (10,0.011077) (14,0.009051) (18,0.00966) (22,0.010341) (26,0.009376) (30,0.007707) (34,0.00809) (38,0.00763) (42,0.006809) (46,0.007939) };
\addplot +[teal] coordinates {(2,0.014737) (6,0.029147) (10,0.002143) (14,0.0) (18,0.0) (22,0.0) (26,0.0) (30,0.0) (34,0.0) (38,0.0) (42,0.0) (46,0.0) };

\legend{RoR,Nuclear,Storage,Biomass,PV,Wind,Others}
\end{axis}
\end{tikzpicture}
	\vspace{-0.3cm}
	\caption{Monthly simulated production per technology in CH}\label{2030 generation}
	\vspace{-0.4cm}
\end{figure}\\
\indent Fig. \ref{2030 generation} compares the past (2015) and future (2030) monthly simulated generation. The largest differences between 2015 and 2030 occur during the winter months (Nov-Mar) when the reactors in 2015 are producing at high levels and there is less PV generation, and in June. In both future scenarios, the remaining Swiss nuclear reactor Leibstadt is shut down for scheduled refueling in June. As a result, the 2015 production can not be reached despite the high solar output and Switzerland becomes a net importer during this month which used to be an export month in 2015. In 2030, the majority of electricity imports come from France as opposed to Germany (in 2015). This is due to the complete nuclear phaseout in Germany and the projected growth of RES in France.



\section{Conclusion and Outlook} \label{C}
This paper presents a GEP formulation which provides the location, type and size of new generators/storages considering system flexibility needs. The multitude of results and their high level of detail (both spatial and temporal) could prove to be useful to TSOs, policy makers and asset owners/operators alike. Future work will focus on including more details in the modeling of the surrounding countries and hydro power in Switzerland. Furthermore, we will investigate a coordinated approach to investments in new generation capacity on the transmission and distribution system levels to improve the modeling of PV integration in Switzerland.

\section{Acknowledgements}
The authors would like to thank the Swiss Federal Office of Energy (SFOE) for supporting the Nexus-e project (Nr. SI/501460).
The views and opinions in this paper are those of the authors and do not necessarily reflect the position of SFOE.

\end{document}